\documentclass[secthm,seceqn,amsthm,ussrhead,12pt]{amsart}
\usepackage{amsmath,latexsym}
\usepackage[english]{babel}
\usepackage[psamsfonts]{amssymb}
\usepackage{times}
\usepackage{cite}
\usepackage{pdflscape} 
\usepackage{ulem}
\usepackage[mathcal]{euscript}
\usepackage{tikz}
\usepackage{hyperref}
\usepackage{cancel}
\usepackage{multirow}
\usetikzlibrary{arrows}
\usepackage{color}

\newtheorem{Th}{Theorem}
\newtheorem{Lem}[Th]{Lemma}

\newtheorem{corollary}[Th]{Corollary}

\renewcommand \a {\mathfrak{a}}

\newcommand \Z {\mathbb{Z}}

\newcommand \GL {\operatorname{GL}}
\newcommand \Aut {\operatorname{Aut}}

\newcommand \Hom {\operatorname{Hom}}

\newcommand \tr {\operatorname{tr}}

\newenvironment{Proof}[1][Proof.]{\begin{trivlist}
\item[\hskip \labelsep {\bfseries #1}]}{\flushright
$\Box$\end{trivlist}}

\begin{document}
	\sloppy

\vspace*{2cm}

{\Large Degenerations of Jordan Superalgebras}

\medskip

\medskip

\medskip

\medskip
\textbf{Mar\'ia Alejandra Alvarez$^{a}$, Isabel Hern\'andez$^{b}$, Ivan Kaygorodov$^{c}$}
\medskip

{\tiny
$^{a}$ Departamento de Matem\'aticas, Universidad de Antofagasta, Chile.

$^{b}$ CONACYT - Centro de Investigaci\'on en Matem\'aticas, A.C.  Unidad M\'erida, M\'exico.

$^{c}$ Universidade Federal do ABC, CMCC, Santo Andr\'{e}, Brazil.
\smallskip

    E-mail addresses:\smallskip

    Mar\'ia Alejandra Alvarez (maria.alvarez@uantof.cl),
    
    Isabel Hern\'andez (isabel@cimat.mx),
    
    Ivan Kaygorodov (kaygorodov.ivan@gmail.com),

}

       \vspace{0.3cm}

\vspace*{2cm}

{\bf Abstract.} 
We describe degenerations of three-dimensional Jordan superalgebras over $\mathbb{C}.$  
In particular, we describe all irreducible components in the corresponding varieties.\smallskip

       \vspace{0.3cm}

       \vspace{0.3cm}

\section{Introduction}

       \vspace{0.3cm}

Contractions of Lie algebras are limiting processes between Lie algebras, which have been studied first in physics \cite{13,10}. For example, classical mechanics is a limiting case of quantum mechanics as $\hbar \to 0,$ described by a contraction of the Heisenberg-Weyl Lie algebras to the Abelian Lie algebra of the same dimension. Description of contractions of low dimensional Lie algebras was given in \cite{nesterpop}. 
The study of contractions and graded contractions of binary algebras has a very big history (see, for example, \cite{c1,c2,c3}).
The study of graded contractions of Jordan algebras and Jordan superalgebras was iniciated in \cite{kp03}.
The first attempt of the study of contractions of $n$-ary algebras stays in the variety of Filippov algebras \cite{deaz}.

In mathematics, often a more general definition of contraction is used, the so-called degeneration. Degenerations are related to deformations. Degenerations of algebras is an interesting subject, which was studied in various papers (see, for example, \cite{CKLO13,BC99,S90,GRH,GRH2,BB09,chouhy,BB14,laur03}). In particular, there are many results concerning degenerations of algebras of low dimensions in a  variety defined by a set of identities. One of the important problems in this direction is the description of the so-called rigid algebras \cite{ikv17}. These algebras are of big interest, since the closures of their orbits under the action of generalized linear group form irreducible components of a variety under consideration (with respect to the Zariski topology). 
There are fewer works in which the full information about degenerations was found for some variety of algebras. This problem was solved 
for two-dimensional pre-Lie algebras in \cite{BB09}, 
for two-dimensional Jordan algebras in \cite{jor2},
for three-dimensional Novikov algebras in \cite{BB14},
for three-dimensional Jordan algebras \cite{gkp17},
for four-dimensional Lie algebras in \cite{BC99}, 
for nilpotent four-dimensional Jordan algebras \cite{contr11},
for nilpotent four-dimensional Leibniz and Zinbiel algebras in \cite{kppv}, 
for nilpotent five- and six-dimensional Lie algebras in \cite{S90,GRH}, 
for nilpotent five- and six-dimensional Malcev algebras in \cite{kpv}, 
and for all $2$-dimensional algebras \cite{kv16}.
On the same time, 
the study of degenerations of superalgebras and graded algebras was iniciated in \cite{deggraa}
for associative case and in \cite{degsulie} for Lie superalgebras.

\section{Definitions and notation}

\subsection{Jordan superalgebras}
Jordan algebras appeared as a tool for studies in quantum mechanic in
the paper of Jordan, von Neumann and Wigner \cite{jnw}.
A commutative algebra is called a {\it  Jordan  algebra}  if it satisfies the identity
$$(x^2y)x=x^2(yx).$$ 
The study of the structure theory and other properties of Jordan algebras was initiated by Albert.
Jordan algebras are related with some questions in 
differentional equations \cite{svi91}, superstring theory \cite{superstring}, 
analysis, operator theory, geometry, mathematical biology, mathematical statistics and physics 
(see, the survey of Iordanescu \cite{radu}).

Let $G$ be the Grassmann algebra over $\mathbb{F}$ given by the generators $1, \xi_1, \ldots , \xi_n, \ldots$  and the defining relations $\xi_i^2=0$ and $\xi_i\xi_j=-\xi_j\xi_i.$ 
The elements $1, \xi_{i_1} \xi_{i_2} \ldots  \xi_{i_k},\ i_1<i_2<\ldots <i_k,$ form a basis of the algebra $G$ over $\mathbb{F}$. Denote by $G_0$ and $G_1$ the subspaces spanned by the products
of even and odd lengths, respectively; then $G$ can be represented as the direct sum of these subspaces,
$G = G_0 \oplus G_1.$ 
Here the relations $G_iG_j \subseteq G_{i+j (mod \ 2)}$, $i,j = 0, 1$, hold. 
In other words, $G$ is a $\mathbb{Z}_2$-graded algebra (or a superalgebra) over $\mathbb{F}.$
Suppose now that $A = A_0 \oplus A_1$ is an arbitrary superalgebra over $\mathbb{F}$. Consider the tensor product $G \otimes A$ of $\mathbb{F}$-algebras. 
The subalgebra
$$G(A) = G_0 \otimes A_0 + G_1  \otimes A_1$$
of $G\otimes A$ is referred to as the Grassmann envelope of the superalgebra $A.$
Let $\Omega$ be a variety of algebras over $\mathbb{F}.$ 
A superalgebra $A = A_0 \oplus A_1$ is referred to as an
$\Omega$-superalgebra if its Grassmann envelope $G(A)$ is an algebra in $\Omega.$ 
In particular, $A = A_0  \oplus  A_1$ is
referred to as a Jordan superalgebra if its Grassmann envelope $G(A)$ is a Jordan algebra. 
The study of Jordan superalgebras has very big history (for example, see \cite{MS,ELS08,MZ01,K10,K12,CK07}).


\subsection{Degenerations}
Given an $(m,n)$-dimensional vector superspace $V=V_0\oplus V_1$, the set 
$$\Hom(V \otimes V,V)=(\Hom(V \otimes V,V))_0\oplus (\Hom(V \otimes V,V))_1$$ is a vector superspace of dimension $m^3+3mn^2$. This space has a structure of the affine variety $\mathbb{C}^{m^3+3mn^2}.$ If we fix a basis $\{e_1,\dots,e_m,f_1,\dots,f_n\}$ of $V$, then any $\mu\in \Hom(V \otimes V,V)$ is determined by $m^3+3mn^2$ structure constants 
$\alpha_{i,j}^k,\beta_{i,j}^q,\gamma_{i,j}^q, \delta_{p,q}^k \in\mathbb{C}$ such that 
$$\mu(e_i\otimes e_j)=\sum\limits_{k=1}^m\alpha_{i,j}^ke_k, 
\quad \mu(e_i\otimes f_p)=\sum\limits_{q=1}^n\beta_{i,p}^qf_q,
\quad \mu(f_p\otimes e_i)=\sum\limits_{q=1}^n\gamma_{p,i}^qf_q,
\quad \mu(f_p\otimes f_q)=\sum\limits_{k=1}^m\delta_{p,q}^ke_k.$$ 
A subset $\mathbb{L}(T)$ of $\Hom(V \otimes V,V)$ is {\it Zariski-closed} if it can be defined by a set of polynomial equations $T$ in the variables 
$\alpha_{i,j}^k,\beta_{i,p}^q, \gamma_{p,i}^q, \delta_{p,q}^k$ ($1\le i,j,k\le m,\ 1\leq p,q\leq n$).

Let $\mathcal{S}^{m,n}$ be the set of all superalgebras of dimension $(m,n)$ defined by the family of polinomial super-identities $T$, understood as a subset $\mathbb{L}(T)$ of an affine variety $\Hom(V\otimes V, V)$. Then one can see that $\mathcal{S}^{m,n}$ is a Zariski-closed subset of the variety $\Hom(V\otimes V, V).$ 
The group $G=(\Aut V)_0\simeq\GL(V_0)\oplus\GL(V_1)$ acts on $\mathcal{S}^{m,n}$ by conjugations:
$$ (g * \mu )(x\otimes y) = g\mu(g^{-1}x\otimes g^{-1}y)$$ 
for $x,y\in V$, $\mu\in\mathbb{L}(T)$ and $g\in G$.

Thus, $\mathcal{S}^{m,n}$ is decomposed into $G$-orbits that correspond to the isomorphism classes of superalgebras. Let $O(\mu)$ denote the orbit of $\mu\in\mathbb{L}(T)$ under the action of $G$ and $\overline{O(\mu)}$ denote the Zariski closure of $O(\mu)$.

Let $J, J' \in \mathcal{S}^{m,n}$  and $\lambda,\mu\in \mathbb{L}(T)$ represent $J$ and $J'$ respectively. We say that $\lambda$ degenerates to $\mu$ and write $\lambda\to \mu$ if $\mu\in\overline{O(\lambda)}$. Note that in this case we have $\overline{O(\mu)}\subset\overline{O(\lambda)}$. Hence, the definition of a degeneration does not depend on the choice of $\mu$ and $\lambda$, and we will right indistinctly $J\to J'$ instead of $\lambda\to\mu$ and $O(J)$ instead of $O(\lambda)$. If $J\not\cong J'$, then the assertion $J\to J'$ is called a {\it proper degeneration}. We write $J\not\to J'$ if $J'\not\in\overline{O(J)}$.

Let $J$ be represented by $\lambda\in\mathbb{L}(T)$. Then $J$ is  {\it rigid} in $\mathbb{L}(T)$ if $O(\lambda)$ is an open subset of $\mathbb{L}(T)$.  Recall that a subset of a variety is called irreducible if it cannot be represented as a union of two non-trivial closed subsets. A maximal irreducible closed subset of a variety is called {\it irreducible component}.  In particular, $J$ is rigid in $\mathcal{S}^{m,n}$ iff $\overline{O(\lambda)}$ is an irreducible component of $\mathbb{L}(T)$. It is well known that any affine variety can be represented as a finite union of its irreducible components in a unique way. We denote by $Rig(\mathcal{S}^{m,n})$ the set of rigid superalgebras in $\mathcal{S}^{m,n}$.

\subsection{Principal notations}
Let $\mathcal{JS}^{m,n}$ be the set of all Jordan superalgebras of dimension $(m,n).$
Let $J$ be a Jordan superalgebra with fixed basis $\{e_1,\dots,e_m,f_1,\dots f_n\}$, defined by
\[e_ie_j=\sum_{k=1}^m\alpha_{ij}^ke_k,\quad e_if_j=\sum_{k=1}^n\beta_{ij}^kf_k,\quad f_if_j=\sum_{k=1}^m\gamma_{ij}^ke_k.\]
We will use the following notation: 
\begin{enumerate}
\item $\a(J)$ is the Jordan superalgebra with the same underlying vector superspace than $J$, and defined by $f_if_j=\displaystyle\sum_{k=1}^n\gamma_{ij}^ke_k$. 
\item $J^1=J$, $J^r=J^{r-1}J+J^{r-2}J^2+\dots+ JJ^{r-1}$, and in every case $J^r=(J^r)_0\oplus (J^r)_1$.
\item $c_{i,j}=\displaystyle\frac{\tr (L(x)^i)\cdot\tr(L(y)^j)}{\tr( L(x)^i\cdot L(y)^j)}$ 
is the Burde invariant, where $L(x)$ 
is the left 
multiplication
. This invariant $c_{i,j}$ 
is defined as a quotient of two polynomials in the structure constants of $J$, for all $x,y\in J$ such that both polynomials are not zero and $c_{i,j}$ 
is independent of the choice of $x,y$.
\end{enumerate}



\section{Methods} 


First of all, if $J\to J'$ and $J\not\cong J'$, then $\dim\Aut(J)<\dim\Aut(J')$, where $\Aut(J)$ is the space of automorphisms of $J$. Secondly, if $J\to J'$ and $J'\to J''$ then $J\to J''$. If there is no $J'$ such that $J\to J'$ and $J'\to J''$ are proper degenerations, then the assertion $J\to J''$ is called a {\it primary degeneration}. If $\dim\Aut(J)<\dim\Aut(J'')$ and there are no $J'$ and $J'''$ such that $J'\to J$, $J''\to J'''$, $J'\not\to J'''$ and one of the assertions $J'\to J$ and $J''\to J'''$ is a proper degeneration,  then the assertion $J \not\to J''$ is called a {\it primary non-degeneration}. It suffices to prove only primary degenerations and non-degenerations to describe degenerations in the variety under consideration. It is easy to see that any superalgebra degenerates to the superalgebra with zero multiplication. From now on we will use this fact without mentioning it.

Let us describe the methods for proving primary non-degenerations. The main tool for this is the following lemma.

\begin{Lem}\label{lema:inv}
If $J\to J'$ then the following hold:
\begin{enumerate}
\item $\dim (J^r)_i\geq\dim (J'^r)_i$, for $i\in\Z_2;$
\item $(J)_0\to (J')_0;$
\item $\a(J)\to\a(J');$
\item If the Burde invariant exist for $J$ and $J'$, then both superalgebras have the same Burde invariant$;$
\item If $J$ is associative then $J'$ must be associative. In fact, if $J$ satisfies a P.I. then $J'$ must satisfy the same P.I.
\end{enumerate}
\end{Lem}

In the cases where all of these criteria can't be applied to prove $J\not\to J'$, we will define $\mathcal{R}$ by a set of polynomial equations and will give a basis of $V$, in which the structure constants of $\lambda$ give a solution to all these equations. We will omit everywhere the verification of the fact that $\mathcal{R}$ is stable under the action of the subgroup of upper triangular matrices and of the fact that $\mu\not\in\mathcal{R}$ for any choice of a basis of $V$. These verifications can be done by direct calculations.

{\bf Degenerations of Graded algebras}.  Let
 $G$ be an abelian group and  let  $\mathcal {V}(\mathcal{F})$ be a variety of algebras defined by a family of  polynomial identities $\mathcal{F}$. It is important to notice that degeneration on the  $G$-graded variety $G\mathcal{V}( \mathcal{F})$ is a more restrictive notion than degeneration on the variety $\mathcal{V}(\mathcal{F})$, In fact, consider $A, A^\prime \in  G\mathcal  {V}(\mathcal{F})$  such that  $ A,  A^\prime  \in \mathcal{V}(\mathcal{F})$, a degeneration between the algebras $A$ and $A^\prime$  may not give rise  to a degeneration  between  the $G$-graded algebras $A$ and $A^\prime$, since  the matrices describing the  basis changes in  $G\mathcal  {V}(\mathcal{F})$  must preserve the $G$-graduation. Hence,  we have the following  result.  

\begin{Lem}
 Let  $ A, A^\prime \in G\mathcal {V}(\mathcal{F}) \cap \mathcal {V}(\mathcal{F})$. If $A  \not \to  A^\prime $ as algebras, then $A \not \to A^\prime $ as $G$-graded algebras. 
\label{alg-nd-Galg-nd}
\end{Lem}

\section{Main result}

In this section we describe all degenerations and non-degenerations of $3$-dimensional Jordan superalgebras.
Note that, 
there is only one (trivial) Jordan superalgebra of the type $(0,3);$
the variety of $3$-dimensional Jordan algebras (Jordan superalgebras of the type $(3,0)$)
has $19$ non isomorphic algebras with non zero multiplication,
in particularly, $5$ algebras are rigid.
The full description of all degenerations and non-degenerations of $3$-dimensional Jordan algebras was given in \cite{gkp17}.
The rest of the section is dedicated to study of degenerations of Jordan superalgebras of types  $(1,2)$ and $(2,1).$

\subsection{Jordan Superalgebras of type $(1,2)$}

\subsubsection{Algebraic classification}
As was noticed the algebraic classification of Jordan superalgebras of the type $(1,2)$ was received in \cite{MS}.
In the next table we give this classification with some additional useful information about these superalgebras.

\begin{center}
Table 1. {\it $(1,2)$-Dimensional Jordan superalgebras.}
\begin{equation*}
\begin{array}{|c|l|c|c|l|} 
\hline
\mbox{$J$}  & \mbox{ multiplication tables } & \mbox{$\dim\Aut(J)$} & \mbox{$c_{i,j}$} & \mbox{type}\\ 
\hline \hline 
U^s_1 & e_1e_1=e_1 & 4 & 1 & \text{associative}\\ \hline
S_1^2 & e_1e_1=e_1, e_1f_1=\frac{1}{2}f_1 & 2 & 2 & \text{non-associative}\\ \hline
S_1^3 & e_1f_1=f_2, \ f_1f_2=e_1                           & 2 & \not\exists & \text{non-associative}\\ \hline
S^2_2 & e_1e_1=e_1, e_1f_1=f_1 & 2 & 2 & \text{associative} \\ \hline
S_2^3 & f_1f_2=e_1            & 4 & \not\exists & \text{associative} \\ \hline
S_3^3 & e_1f_1=f_2                         & 3 & \not\exists & \text{associative}\\ \hline
S^3_4 & e_1e_1=e_1, e_1f_1=f_1, e_1f_2=\frac{1}{2}f_2 & 2 & \frac{\left(2+\left(\frac{1}{2}\right)^i\right)\left(2+\left(\frac{1}{2}\right)^j\right)}{\left(2+\left(\frac{1}{2}\right)^{i+j}\right)} & \text{non-associative}\\ \hline
S^3_5 & e_1e_1=e_1, e_1f_1=\frac{1}{2}f_1, e_1f_2=\frac{1}{2}f_2 & 4 & \frac{\left(1+2\left(\frac{1}{2}\right)^i\right)\left(1+2\left(\frac{1}{2}\right)^j\right)}{\left(1+2\left(\frac{1}{2}\right)^{i+j}\right)} & \text{non-associative}\\ \hline
S^3_6 & e_1e_1=e_1, e_1f_1=f_1, e_1f_2=f_2 & 4 & 3 & \text{associative}\\ \hline
S^3_7 & e_1e_1=e_1, e_1f_1=\frac{1}{2}f_1, e_1f_2=\frac{1}{2}f_2, f_1f_2=e_1 & 3 & \frac{\left(1+2\left(\frac{1}{2}\right)^i\right)\left(1+2\left(\frac{1}{2}\right)^j\right)}{\left(1+2\left(\frac{1}{2}\right)^{i+j}\right)} & \text{non-associative} \\ \hline
S^3_8 & e_1e_1=e_1, e_1f_1=f_1, e_1f_2=f_2, f_1f_2=e_1 & 3 & 3 & \text{non-associative} \\ \hline

\end{array}
\end{equation*}
\end{center}


%

\subsubsection{Degenerations}

\begin{Th}\label{third}\label{theorem}
The graph of primary degenerations for Jordan superalgebras of dimension $(1,2)$ has the following form: 
\end{Th}

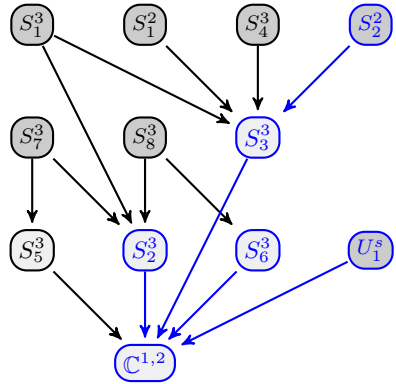
\scriptsize
 
\begin{center}

\begin{tikzpicture}[->,>=stealth',shorten >=0.08cm,auto,node distance=1.5cm,
                    thick,main node/.style={rectangle,draw,fill=gray!12,rounded corners=1.5ex,font=\sffamily \bf \bfseries },
                    blue node/.style={rectangle,draw, color=blue,fill=gray!12,rounded corners=1.5ex,font=\sffamily \bf \bfseries },
                    orange node/.style={rectangle,draw, color=orange,fill=gray!12,rounded corners=1.5ex,font=\sffamily \bf \bfseries },
                    green node/.style={rectangle,draw, color=green,fill=gray!12,rounded corners=1.5ex,font=\sffamily \bf \bfseries },
                    olive node/.style={rectangle,draw, color=olive,fill=gray!12,rounded corners=1.5ex,font=\sffamily \bf \bfseries },
                    connecting node/.style={circle, draw, color=purple },
                    rigid node/.style={rectangle,draw,fill=black!20,rounded corners=1.5ex,font=\sffamily \tiny \bfseries },
                    bluerigid node/.style={rectangle,draw,color=blue,fill=black!20,rounded corners=1.5ex,font=\sffamily \tiny \bfseries                    }, style={draw,font=\sffamily \scriptsize \bfseries }]
                    
\node (30)   {};

\node (31)[right of=30]{};	
\node (32)[right of=31]{};	
\node (33)[right of=32]{};	
\node (34)[right of=33]{};

\node  [rigid node] (s31) [right of=30] {$S^3_1$};	
\node  [rigid node] (s21s11) [right of=31] {$S^2_1$};	
\node  [rigid node] (s34) [right of=32] {$S^3_4$};	
\node  [bluerigid node] (s22s11) [right of=33] {$S^2_2$};	

\node (20)[below of=30]{};	
\node (21)[right of=20]{};	
\node (22)[right of=21]{};	
\node (23)[right of=22]{};	
\node (24)[right of=23]{};	

\node  [blue node] (s33) [right of=22] {$S^3_3$};	
\node  [rigid node] (s37) [right of=20] {$S^3_7$};	
\node  [rigid node] (s38) [right of=21] {$S^3_8$};

\node (10)[below of=20]{};	
\node (11)[right of=10]{};	
\node (12)[right of=11]{};	
\node (13)[right of=12]{};	
\node (14)[right of=13]{};	

\node  [blue node] (s32) [right of=11] {$S^3_2$};	
\node  [bluerigid node] (u1ss) [right of=13] {$U_1^s$};	
\node  [main node] (s35) [right of=10] {$S^3_5$};	
\node  [blue node] (s36) [right of=12] {$S^3_6$};

\node (00)[below of=10]{};	
\node (01)[right of=00]{};	
\node (02)[right of=01]{};	

\node  [blue node] (u2ss) [right of=01] {$\mathbb{C}^{1,2}$};

\path[every node/.style={font=\sffamily\small}]


(s31) edge [bend right=0, color=black] node{}  (s33)
(s31) edge [bend right=0, color=black] node{}  (s32)

(s21s11) edge [bend right=0, color=black] node{}  (s33)

(s34) edge [bend right=0, color=black] node{}  (s33)
(s22s11) edge [bend right=0, color=blue] node{}  (s33)
(s37) edge [bend right=0, color=black] node{}  (s32)
(s37) edge [bend right=0, color=black] node{}  (s35)
(s38) edge [bend right=0, color=black] node{}  (s32)
(s38) edge [bend right=0, color=black] node{}  (s36)
(s32) edge [bend right=0, color=blue] node{}  (u2ss)
(u1ss) edge [bend right=0, color=blue] node{}  (u2ss)
(s35) edge [bend right=0, color=black] node{}  (u2ss)
(s36) edge [bend right=0, color=blue] node{}  (u2ss)
(s33) edge [bend right=0, color=blue] node{}  (u2ss);

\end{tikzpicture}

\end{center}


\normalsize
 

\begin{Proof} 

We prove all required primary degenerations in Table 2 below. Recall that an associative superalgebra can only degenerate to an associative superalgebra. Let us consider the first degeneration $S^2_1\to S_3^3$ 
to clarify this table. Write nonzero products in $S^2_1$ in the basis  $E_1^t,F_1^t,F_2^t$: 
$$E_1^tE_1^t=tE_1^t, \ E_1^tF_1^t=te_1(f_1-2t^{-1}f_2)=\frac{t}{2}(f_1-2t^{-1}f_2)+f_2=\frac{t}{2}F_1^t+F_2^t.$$
 It is easy to see that for $t=0$ we obtain the multiplication table of $S_3^3.$ 
 The remaining degenerations can be considered in the same way.

\begin{center}\footnotesize Table 2. {\it Primary degenerations of Jordan superalgebras of dimension $(1,2)$.}
$$\begin{array}{|c|lll|}
\hline
\mbox{degenerations}  &  \multicolumn{3}{|c|}{\mbox{parametrized bases}} \\
\hline
\hline
S^2_1\to S_3^3   & E^t_1=te_1,& F_1^t=f_1-2t^{-1}f_2,& F_2^t=f_2  \\ \hline
S^2_2\to S_3^3   & E^t_1=te_1,& F_1^t=f_1-t^{-1}f_2,& F_2^t=f_2  \\ \hline
S^3_1\to S^3_3   & E^t_1=e_1,& F_1^t=tf_1,& F_2^t=tf_2  \\ \hline
S^3_1\to S^3_2   & E^t_1=e_1,& F_1^t=tf_1,& F_2^t=tf_2  \\ \hline
S^3_4\to S^3_3   & E^t_1=te_1,& F_1^t=f_1-2t^{-1}f_2,& F_2^t=f_2  \\ \hline
S^3_7\to S^3_2   & E^t_1=te_1, & F_1^t=tf_1,& F_2^t=f_2  \\ \hline
S^3_7\to S^3_5   & E^t_1=e_1, &F_1^t=f_1,& F_2^t=tf_2  \\ \hline
S^3_8\to S^3_2   & E^t_1=te_1,& F_1^t=tf_1,& F_2^t=f_2  \\ \hline
S^3_8\to S^3_6   & E^t_1=e_1,& F_1^t=tf_1,& F_2^t=f_2  \\ \hline
\end{array}$$

\end{center}

\normalsize

The primary non-degenerations are proved in Table 3.

\begin{center}\footnotesize Table 3. {\it Primary non-degenerations of Jordan superalgebras of dimension $(1,2)$.}
$$\begin{array}{|l|c|}
\hline
\mbox{non-degenerations 
}  &  \mbox{reason} \\
\hline
\hline
\begin{array}{l}
S^3_3\not\to S_2^3 \end{array}& \dim (J^2)_0<\dim (J'^2)_0 \\ \hline

\begin{array}{l}
S^3_3\not\to U^s_1,\ 
S^3_6;\quad 
S^3_1\not\to 
S^3_7,\ 
S^3_8,\ 
U^s_1,\ 
S^3_5,\ 
S^3_6 \end{array}
& J_0\not\to J'_0 \\ \hline

\begin{array}{l}

S^3_7\not\to U_1^s,\ 
S_6^3;\quad 
S^3_8\not\to U_1^s,\ 
S_3^5;\quad  

S^2_1\not\to S_7^3,\ 
S_8^3,\ 
U_1^s,\ 
S_5^3,\ 
S_6^3; \\

S^3_4\not\to S^3_7,\ 
S^3_8,\ 
U^s_1,\ 
S^3_5,\ 
S^3_6;\quad 

S^2_2\not\to U_1^s,\ 
S_6^3 
\end{array}
&   c_{i,j} 
\\ \hline

\begin{array}{l}
S^2_1\not\to S_2^3;\quad 
S^3_4\not\to S_2^3;\quad
S^2_2\not\to S_2^3;\quad
\end{array}& \a(J)\not\to\a(J')  \\ \hline
\end{array}$$

\end{center}

 \end{Proof}

\subsubsection{Irreducible components and rigid algebras}

Using Theorem \ref{theorem}, we describe the irreducible components and the rigid algebras in $\mathcal{JS}^{1,2}.$

\begin{corollary}\label{ir_12} The irreducible components of $\mathcal{JS}^{1,2}$ are:

$$
\begin{aligned}
\mathcal{C}_1   &=\overline{ O(U_{1}^s)  }=  
\{  U_1^s,  \mathbb{C}^{1,2} \}:\\  
\mathcal{C}_2   &=\overline{ O(S_{1}^2)  }=  
\{  S_1^2, S_3^3, \mathbb{C}^{1,2} \}; \\
\mathcal{C}_3   &=\overline{ O(S_{1}^3)  }=  
\{  S_1^3, S_3^3, \mathbb{C}^{1,2} \}; \\
\mathcal{C}_4   &=\overline{ O(S_{2}^2)  }=  
\{  S_2^2, S_3^3, S^3_2, \mathbb{C}^{1,2} \}; \\
\mathcal{C}_5   &=\overline{ O(S_{4}^3)  }=  
\{  S_4^3, S^3_3, \mathbb{C}^{1,2} \}; \\
\mathcal{C}_6   &=\overline{ O(S_{7}^3)  }=  
\{  S_5^3, S_2^3, \mathbb{C}^{1,2} \}; \\
\mathcal{C}_7   &=\overline{ O(S_{8}^3)  }=  
\{  S_8^3, S_2^3, S^3_6, \mathbb{C}^{1,2} \};\\
\end{aligned}
$$ 

In particular, $Rig(\mathcal{JS}^{1,2})= \{ U_{1}^s, S_{1}^2, S_{1}^3, S_{2}^2, S_{4}^3, S_{7}^3, S_{8}^3 \}.$

\end{corollary}

\subsection{Superalgebras of type $(2,1)$}

In this section we describe all possible primary degenerations between Jordan superalgebras of dimension $(2,1)$. 
First of all, notice that every Jordan superalgebra of dimension $(2,1)$ has trivial odd product, so it can be considered as a Jordan algebra of dimension $3$. However, the graph of primary degenerations of Jordan superalgebras of dimension $(2,1)$ is not a subgraph of the primary degenerations of de Jordan algebras of dimension $3$. In fact,  there exist Jordan superalgebras, without degenerations
between them,  such that they degenerate as Jordan algebras. Take for example
$S_2^2 \oplus U_1^s$ and $B_1^s$ (see Table 5),  Taking the basis change:  $E_1^t= e_1$, 
$E_2^t=f_1$ and $F_1^t= t e_2$, it follows that $S_2^2 \oplus U_1^s  \to B_1^s$ as Jordan algebras. 
Notice that this basis change does not preserve the $\mathbb{Z}_2$-graduation. Moreover,  we shall prove that 
$S_2^2 \oplus U_1^s \not \to B_1^s$ as Jordan superalgebras.


\subsubsection{Algebraic classification}
In the next table we provide the classification of $(2,1)$-dimensional Jordan superalgebras with some additional useful information about these superalgebras.

\begin{center}
Table 4. {\it $(2,1)$-Dimensional Jordan superalgebras.}

\begin{equation*}
\begin{array}{|l|l|c|l|} 
\hline
\mbox{$J$}  & \mbox{Multiplication tables } & \mbox{$\dim\Aut(J)$} & \mbox{Type} \\ 
\hline \hline 
2U_1^s 
& e_1e_1=e_1, \; \;e_2e_2=e_2 & 1 & \mbox{associative} \\\hline

U_1^s 
& e_1e_1=e_1 & 2 & \mbox{associative} \\ \hline

B_1^s 
& e_1e_1=e_1, \;\; e_1e_2=e_2 & 2 & \mbox{associative} \\ \hline
B_2^s 
& e_1e_1=e_1, \;\;e_1e_2= \frac{1}{2}e_2 & 3 & \mbox{non-associative} \\ \hline

B_3^s 
& e_1e_1=e_2   &  3 & \mbox{associative} \\ \hline

S_1^2\oplus U_1^s & e_1e_1=e_1, \;\; e_2e_2=e_2, \;\; e_1f_1=\frac{1}{2}f_1 & 1 & \mbox{non-associative} \\ \hline

S_1^2 
& e_1e_1=e_1,\;\; e_1f_1= \frac{1}{2}f_1 & 2 & \mbox{non-associative} \\ \hline

S_2^2 \oplus U_1^s & e_1e_1=e_1, \;\;e_2e_2=e_2,  \;\;e_1f_1=f_1 & 1  & \mbox{associative} \\ \hline

S_2^2 
& e_1e_1=e_1, \;\; e_1f_1=f_1 & 2 &  \mbox{associative}\\ \hline
S_9^3  & e_1e_1=e_1,\;\; e_1e_2=e_2,  \;\; e_1f_1=\frac{1}{2}f_1 & 2 & \mbox{non-associative} \\ \hline
S_{10}^3  & e_1e_1=e_1, \;\; e_1e_2=e_2, \;\;e_1f_1=f_1 & 2 & \mbox{associative} \\ \hline
S_{11}^3   & e_1e_1=e_1, \;\;e_1e_2= \frac{1}{2}e_2, \;\;e_1f_1= \frac{1}{2}f_1 & 3 & \mbox{non-associative} \\ \hline
S_{12}^3   & e_1e_1=e_1,\;\; e_1e_2= \frac{1}{2}e_2, \; e_1f_1=f_1 & 3 & \mbox{non-associative} \\ \hline
S_{13}^3 & e_1e_1=e_1, \;\;e_2e_2=e_2, \;\; e_1f_1=\frac{1}{2}f_1,\;\; e_2f_1=\frac{1}{2}f_1 & 1 & \mbox{non-associative}\\ \hline
\end{array}
\label{2-1JSA}
\end{equation*}
\end{center}

\subsubsection{Degenerations}

\begin{Th}\label{4th}\label{theorem21}
The graph of primary degenerations for Jordan superalgebras of dimension $(2,1)$ has the following form: 
\end{Th}

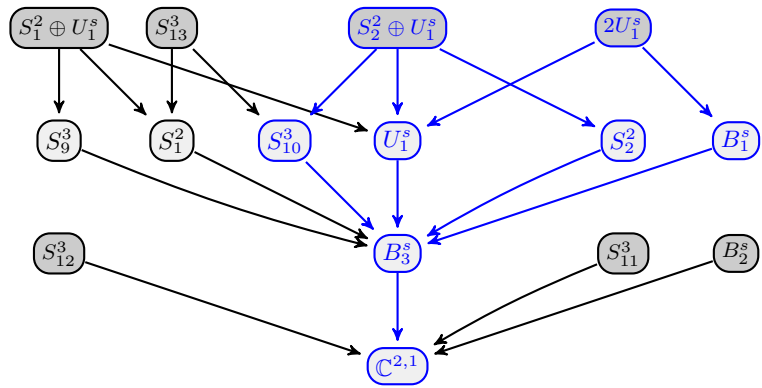
\begin{center}
\scriptsize
\begin{tikzpicture}[->,>=stealth',shorten >=0.08cm,auto,node distance=1.5cm,
                    thick,main node/.style={rectangle,draw,fill=gray!12,rounded corners=1.5ex,font=\sffamily \bf \bfseries },
                    blue node/.style={rectangle,draw, color=blue,fill=gray!12,rounded corners=1.5ex,font=\sffamily \bf \bfseries },
                    orange node/.style={rectangle,draw, color=orange,fill=gray!12,rounded corners=1.5ex,font=\sffamily \bf \bfseries },
                    green node/.style={rectangle,draw, color=green,fill=gray!12,rounded corners=1.5ex,font=\sffamily \bf \bfseries },
                    olive node/.style={rectangle,draw, color=olive,fill=gray!12,rounded corners=1.5ex,font=\sffamily \bf \bfseries },
                    connecting node/.style={circle, draw, color=purple },
                    rigid node/.style={rectangle,draw,fill=black!20,rounded corners=1.5ex,font=\sffamily \tiny \bfseries },
                    bluerigid node/.style={rectangle,draw,color=blue,fill=black!20,rounded corners=1.5ex,font=\sffamily \tiny \bfseries                    },
                    style={draw,font=\sffamily \scriptsize \bfseries }]

\node (30)   {};

\node (31)[right of=30]{};	
\node (32)[right of=31]{};	
\node (33)[right of=32]{};	
\node (34)[right of=33]{};	
\node (35)[right of=34]{};	
\node (36)[right of=35]{};

\node  [rigid node] (s12u1s) [right of=30] {$S_1^2 \oplus U_1^s$};	
\node  [rigid node] (s133) [right of=31] {$S_{13}^3$};	
\node  [bluerigid node] (s22u1s) [right of=33] {$S_2^2 \oplus U_1^s$};	
\node  [bluerigid node] (u1su1ss11) [right of=35] {$2U_1^s 
$};

\node (20)[below of=30]{};	
\node (21)[right of=20]{};	
\node (22)[right of=21]{};	
\node (23)[right of=22]{};	
\node (24)[right of=23]{};	
\node (25)[right of=24]{};	
\node (26)[right of=25]{};	
\node (27)[right of=26]{};	
\node (28)[right of=27]{};

\node  [main node] (s93) [right of=20] {$S^3_9$};	
\node  [main node] (s12u2s) [right of=21] {$S_1^2 
$};
\node  [blue node] (s103) [right of=22] {$S^3_{10}$};
\node  [blue node] (u1su2ss11) [right of=23] {$U_1^s 
$};
\node  [blue node] (s22u2s) [right of=25] {$S_2^2 
$};	
\node  [blue node] (b1ss12) [right of=26] {$B_1^s 
$};

\node (10)[below of=20]{};	
\node (11)[right of=10]{};	
\node (12)[right of=11]{};	
\node (13)[right of=12]{};	
\node (14)[right of=13]{};	
\node (15)[right of=14]{};	
\node (16)[right of=15]{};

\node  [rigid node] (s123) [right of=10] {$S_{12}^3$};	
\node  [blue node] (b3ss11) [right of=13] {$B_3^s 
$};	
\node  [rigid node] (b2ss11) [right of=16] {$B_2^s 
$};	
\node  [rigid node] (s113) [right of=15] {$S_{11}^3$};	

\node (00)[below of=10]{};	
\node (01)[right of=00]{};	
\node (02)[right of=01]{};	
\node (03)[right of=02]{};	
\node (04)[right of=03]{};

\node  [blue node] (u2su2ss11) [right of=03] {$\mathbb{C}^{2,1}     $};

\path[every node/.style={font=\sffamily\small}]


(s133) edge [bend right=0, color=black] node{}  (s12u2s)
(s133) edge [bend right=0, color=black] node{}  (s103)
(s12u1s) edge [bend right=0, color=black] node{}  (s12u2s)
(s12u1s) edge [bend right=0, color=black] node{}  (s93)
(s12u1s) edge [bend right=0, color=black] node{}  (u1su2ss11)
(s12u2s) edge [bend right=0, color=black] node{}  (b3ss11)
(s93) edge [bend right=4, color=black] node{}  (b3ss11)
(b2ss11) edge [bend right=0, color=black] node{}  (u2su2ss11)
(s113) edge [bend right=4, color=black] node{}  (u2su2ss11)
(s123) edge [bend right=0, color=black] node{}  (u2su2ss11)
(u1su1ss11) edge [bend right=0, color=blue] node{}  (b1ss12)
(u1su1ss11) edge [bend right=0, color=blue] node{}  (u1su2ss11)
(s22u1s) edge [bend right=0, color=blue] node{}  (u1su2ss11)
(s22u1s) edge [bend right=0, color=blue] node{}  (s22u2s)
(s22u1s) edge [bend right=0, color=blue] node{}  (s103)
(b1ss12) edge [bend right=0, color=blue] node{}  (b3ss11)
(u1su2ss11) edge [bend right=0, color=blue] node{}  (b3ss11)
(s22u2s) edge [bend right=4, color=blue] node{}  (b3ss11)
(s103) edge [bend right=0, color=blue] node{}  (b3ss11)
(b3ss11) edge [bend right=0, color=blue] node{}  (u2su2ss11);

\end{tikzpicture}
\end{center}

\begin{Proof} We prove all required primary degenerations in Table 5 below. 
For clarifying this table was conside\-red an example in the proof of the theorem \ref{theorem}.

\begin{center}\footnotesize Table 5. {\it Primary degenerations of Jordan superalgebras of dimension $(2,1)$.}
$$\begin{array}{|l|lll|}
\hline
\mbox{degenerations}  &  \multicolumn{3}{|c|}{\mbox{parametrized bases}} \\
\hline
\hline
2U_1^s 
\to  U_1^s 
& E^t_1=e_1,& E_2^t=te_2,& F_1^t=f_1\\ \hline

2U_1^s 
\to B_1^s 
& E^t_1=e_1+e_2,& E_2^t=te_2,& F_1^t=f_1\\ \hline

U_1^s 
\to  B_3^s 
& E^t_1=te_1+e_2,& E_2^t=t^2e_1,& F_1^t=f_1\\ \hline

B_1^s 
\to  B_3^s 
& E^t_1=te_1+e_2,& E_2^t=-t^2e_1,& F_1^t=f_1\\ \hline

S_1^2  \oplus U_1^s \to U_1^s 
& E^t_1=e_2,& E_2^t=te_1,& F_1^t=f_1\\ \hline

S_1^2  \oplus U_1^s \to S_1^2 
& E^t_1=e_1,& E_2^t=te_2,& F_1^t=f_1\\ \hline

S_1^2  \oplus U_1^s \to S_9^3 & E^t_1=e_1+e_2,& E_2^t=te_2,& F_1^t=f_1\\ \hline

S_1^2  
\to  B_3^s 
& E^t_1=te_1+e_2,& E_2^t=-t^2e_1,& F_1^t=f_1\\ \hline

S_2^2  
\to  B_3^s 
& E^t_1=te_1+e_2,& E_2^t=t^2e_1,& F_1^t=f_1\\ \hline

S_2^2  \oplus U_1^s \to S_2^2 
& E^t_1=e_1,& E_2^t=te_2,& F_1^t=f_1\\ \hline

S_2^2  \oplus U_1^s \to U_1^s 
& E^t_1= \frac{1}{t}e_2,& E_2^t=e_1,& F_1^t=f_1\\ \hline

S_2^2  \oplus U_1^s \to S_{10}^3 & E^t_1=e_1+e_2,& E_2^t=te_2,& F_1^t=f_1\\ \hline

S_{9}^3 \to  B_3^s 
& E^t_1=te_1+e_2,& E_2^t=te_2, & F_1^t=f_1\\ \hline

S_{10}^3 \to  B_3^s 
& E^t_1=te_1-\frac{1}{2}e_2,& E_2^t=t^2e_1- te_2   ,& F_1^t=f_1\\ \hline

S^3_{13}\to S_1^2 
& E^t_1=e_1,& E_2^t=te_2,& F_1^t=f_1  \\ \hline

S^3_{13}\to S_{10}^3 & E^t_1=e_1+e_2,& E_2^t=te_2,& F_1^t=f_1  \\ \hline

\end{array}$$
\end{center}

Primary non-degenerations between 3-dimensional Jordan algebras are given in \cite{gkp17}, we use this result and Lemma 2 to prove some  primary non-degenerations, (see table $6$).   

From \cite{gkp17} it follows that  $2U_1^s  \to S_2^2 $, $ S_1^2 \to B_2^s$, $S_9^3 \to S_{12}^3$, and   $S_2^2 \oplus U_1^s \to B_1^s$ as Jordan algebras; we shall prove that they do not degenerate as Jordan superalgebras.

First,  suppose that  $S_2^2 \oplus U_1^s \to B_1^s$  as Jordan superalgebras,  then there exists a parameterized basis 
$$E_1^t= a(t)e_1+ b(t)e_2,   \quad E_2^t= c(t)e_1 + d(t) e_2,   \quad F_1^t= x(t)f_1$$ 
for $S_2^2 \oplus U_1^2$, such that for $t=0$ we obtain  $B_1^s$. 
Since  $E_1^t F_1^t = a(t) F_1^t$ and $E_2^t F_1^t = c(t)  F_1^t$, it follows that $a(0)=c(0)=0$.  Now, since 
$E_1^tE_1^t=E_1^t$ it follows that  $a(t)=0$ and $b(t)=1$ for all $t$,  Finally, since $E_2^tE_2^t =0$ at $t=0$, it follows that $d(0)=0$, showing that  $ E_1^t E^t_2 = 0$,  for all $t$.  This proves that $S_2^2 \oplus U_1^s \not \to B_1^s$.  For the remaining cases see table $6$. 


\begin{center}\footnotesize Table 6. {\it Primary non-degenerations of Jordan superalgebras of dimension $(2,1)$.}
$$\begin{array}{|l|c|}
\hline
\mbox{non-degenerations 
}  &  \mbox{reason} \\
\hline
\hline

\begin{array}{l}
2U_1^s 
\not \to S_1^2, \ S_2^2, \ S_9^3, \ S_{10}^3, \  S_{11}^3, \ S_{12}^3    
\end{array}
 & \dim (J^r)_1<\dim ((J^\prime)^r)_1   \\ \hline 
 
\begin{array}{l}
S_1^2 
 \not \to B_2^s 
;\quad S_9^3 \not \to  S_{12}^3;\quad S_1^2 \oplus U_1^s\not \to B_2^s, \ S_{11}^3, \ S_{12}^3; \\
S_{13}^3 \not \to B_2^s, \ S_{11}^3, \ S_{12}^3 \\
\end{array}      

&  J_0\not\to J'_0 
\\ \hline

\begin{array}{l}
2U_1^s \not \to B_{2}^s; 
S_1^2 \not \to S_{11}^3,\;  S_{12}^3; \quad S_9^3 \not \to B_2^s, \;  S_{11}^3;\\
S_1^2\oplus U_1^s \not \to  B_1^s, \;  S_2^2, \; S_{10}^3; \quad 
S_2^2\oplus U_1^s \not \to  S_{1}^2, \ S_{9}^3; \
S_{13}^3 \not\to U_1^s,\;  B_1^s, \; S_2^2, \; S_9^3
\end{array}
 & J\not \to J' \text{ as Jordan algebras}
\\ \hline
\end{array}$$
\end{center}

\scriptsize

\normalsize

 \end{Proof}

\subsubsection{Irreducible components and rigid algebras}

Using Theorem \ref{theorem21}, we describe the irreducible components and the rigid superalgebras in $\mathcal{JS}^{2,1}.$

\begin{corollary}\label{ir_LZ} The irreducible components of $\mathcal{JS}^{2,1}$ are:

$$
\begin{aligned}
\mathcal{C}_1   &=\overline{ O(2U_1^s)  }=  
\{ 2U_1^s ,  U_1^s , B_1^s, B_3^s,  \mathbb{C}^{2,1} \};\\  
\mathcal{C}_2   &=\overline{ O(B_2^s)  }=  \{ B_2^s ,  \mathbb{C}^{2,1} \};\\  
\mathcal{C}_3   &=\overline{ O(S_{1}^2 \oplus U_1^s)  }=  
\{U_1^s,  B_3^s, S_{1}^2 \oplus U_1^s,   S_1^2 ,   S_9^3, \mathbb{C}^{2,1} \};\\  
\mathcal{C}_4   &=\overline{ O(S_{2}^2 \oplus U_1^s)  }=  
\{  U_1^s,  B_3^s, S_{2}^2 \oplus U_1^s,  S_2^2 , S_{10}^3,  \mathbb{C}^{2,1} \};\\  
\mathcal{C}_5   &=\overline{ O(S_{11}^3)  }=  \{ S_{11}^3,  \mathbb{C}^{2,1} \};\\  
\mathcal{C}_6   &=\overline{ O(S_{12}^3)  }=  \{ S_{12}^3,  \mathbb{C}^{2,1} \};\\  
\mathcal{C}_7   &=\overline{ O(S_{13}^3)  }=  
\{ B_3^s , S_{1}^2 ,  S_{10}^3,    S_{13}^3,  \mathbb{C}^{2,1} \}.\\  
\end{aligned}
$$ 

In particular, $Rig(\mathcal{JS}^{2,1})= 
\{ 2U_1^s, B_2^s , S_{1}^2 \oplus U_1^s,
S_{2}^2 \oplus U_1^s, S_{11}^3, S_{12}^3, S_{13}^3 \}.$

\end{corollary}

\paragraph{Acknowledgements} 
This work was started during the research stay of I. Kaygorodov at the Department of Mathematics, parcially funded by the {\it Coloquio de Matem\'atica} (CR 4430) of the University of Antofagasta. 
The work was supported by RFBR 17-01-00258, by FOMIX-CONACYT  YUC-2013-C14-221183 and 222870.

\end{document}